\documentclass[12pt]{article}
\usepackage[centertags]{amsmath}
\usepackage{amsmath}
\usepackage{amsthm}
\usepackage{color} 
\usepackage[english]{babel}
\usepackage[latin1]{inputenc}
\usepackage[T1]{fontenc}
\usepackage{graphicx} 
\usepackage{amsfonts}
\usepackage{amssymb}
\usepackage[all,cmtip]{xy}

\newtheorem{theorem}{Theorem}[section]

\newtheorem{lemma}[theorem]{Lemma}

\theoremstyle{definition}
\newtheorem{definition}[theorem]{Definition}
\theoremstyle{remarque}

\theoremstyle{example}

\numberwithin{equation}{section}

\begin{document}
\begin{center}\Large\textbf{Completely solvable Lie foliations}
\end{center}
\begin{center}
 {\bf  Ameth \, Ndiaye}\footnote{\tiny\bf{
 Universit\'e\, Cheikh \,Anta \,Diop, \,Dakar/\,D\'epartement de Math\'ematiques(FASTEF)\\Email:\,ameth1.ndiaye@ucad.edu.sn\\
 }}
\end{center}

{\bf Keywords:} foliation, solvable, Lie group, manifold, compact.
\\

\begin{abstract}
In this paper we try to generalize the Haefliger theorem on completly solvable Lie foliations. We prove that: every completely solvable Lie foliation on a compact manifold is the inverse image of a homogene foliation.\\
Every manifold in this paper is compact and our Lie group $G$ is connexe and simply connexe.
\end{abstract}

\section{Introduction}
A codimension $n$ \textbf{foliation} $\mathcal{F}$ on a $(n+m)$-manifold $M$ is given by an open cover $\{U_i\}_{i\in I}$ and submersions
$f_i:U_i\longrightarrow T$ over an $n$-dimensional manifold $T$ and, for $U_i\cap U_j\neq\emptyset$, a diffeomorphism $\gamma_{ij}:f_i(U_i\cap U_j)\longrightarrow f_j(U_i\cap U_j)$ such that $f_j=\gamma_{ij}\circ f_i$.\\
We say that $\{U_i,f_i,T,\gamma_{ij}\}$ is a foliated cocycle defining $\mathcal{F}$.\\
A \textbf{transverse structure} to $\mathcal{F}$ is a geometric structure on $T$ invariant by the local diffeomorphisms $\gamma_{ij}$. We say that $\mathcal{F}$ is a \textbf{Lie $G$-foliation}, if $T$ is a Lie group $G$ and $\gamma_{ij}$ are restrictions of left translations on $G$.\\
Such foliation can also be defined by a $1$-form $\omega$ on $M$ with values in the Lie algebra $\mathcal{G}$ such that:
\begin{enumerate}
\item $\omega_x:T_xM\longrightarrow\mathcal{G}$ is surjective for every $x\in M$
\item  $d\omega+\frac{1}{2}[\omega,\omega]=0$.
\end{enumerate}
In the general case, the structure of a Lie foliation on a compact manifold, is given by the following theorem due to E. F\'edida \cite{2}:\\
\begin{theorem}{\cite{2}}
Let $\mathcal{F}$ be a Lie $G$-foliation on a compact manifold $M$. Let $\widetilde{M}$ be the universal covering of $M$ and  $\widetilde{\mathcal{F}}$ the lift of $\mathcal{F}$ to $\widetilde{M}$. Then there exist a homomorphism $h : \pi_1(M)\longrightarrow G$ and a locally trivial fibration $D:\widetilde{M}\longrightarrow G$ whose fibres are the leaves of $\widetilde{\mathcal{F}}$ and such that, for every $\gamma\in\pi_1(M)$, the following
diagram is commutative:
\end{theorem}
where the first line denotes the deck transformation of $\gamma\in\pi_1(M)$ on $\widetilde{M}$.\\
The group $\Gamma=h(\pi_1(M))$ (which is a subgroup of $G$) is called the holonomy group of $\mathcal{F}$ although the holonomy of each leaf is trivial. The fibration : $D:\widetilde{M}\longrightarrow G$ is called the developing map of $\mathcal{F}$.
The topology of a Lie foliation of a compact manifold is trivial since the leaves are all diffeomorphic, see \cite{2}, which is why the theory of Lie foliation has developed mainly in the sense of classification. This paper is a contribution in this direction.\\
An important class of Lie foliation is the Lie foliation, which is called homogeneous Lie foliation, obtained in the following way: $G$ and $H$ two connected, simply connected, Lie groups. $\Gamma$ a cocompact discrete subgroup of $H$, and $\varphi$ a surjective morphism of Lie groups of $H$ in $G$; the classes on the left of $H$ modulo $\ker\varphi$ reprojected on $G/\Gamma$ are the leaves of a $G$-foliation of $H/\Gamma$ whose holonomy morphism is the restriction of $\varphi$ to $\Gamma$ and the developing application is $\varphi$. We would like to know in which condition a Lie foliation can be reduced, using naturals operations, to a homogeneous Lie foliation. Since the holonomy group of a $G$-foliation reflects the transverse structure, this question is closely related to the following: Given a Lie group $G$, which subgroups are feasible as holonomy group of a $G$-foliation. When $G$ is nilpotent Haefliger showed that in a nilpotent and simply connected $G$ Lie group, any $\Gamma$ subgroup of finite type is a holonomy group of a homogeneous $G$-foliation on a compact manifold. In particular, nilpotent Lie foliations on compact manifolds are all, in reverse order, homogeneous Lie foliations. In the case where the Lie group is resolvable (not necessarily nilpotent), Gaël Meigniez showed that if $G$ a resolvable, connected and simply connected Lie group. $\Gamma$ a subgroup of $G$ of finite type. If $\Gamma$ contains a subgroup that is both polycyclic and uniform, then $\Gamma$ is a holonomy group of $G$-foliation. Our contribution in this paper is to be able to generalise  Haefliger's question in the case of completely solvable Lie groups.
\section{Completely solvable Lie groups}
\begin{definition}
A solvable Lie group $G$ of Lie algebra $\mathcal{G}$ is said to be completely solvable if all the adjoint linear operators $ad_X$ $(X\in\mathcal{G})$ only eigen rea values.
\end{definition}
For example, any nilpotent Lie group is completely solvable.
\begin{definition}
A Lie group $G$ is said to be polycyclic if it has a sequence of cyclic quotient composition. If $G$ is solvable and is of finite type with $\gamma_1, ..., \gamma_n$ as the generating part, then $G$ is polycyclic if and only if the roots of $\gamma_i$ are algebraic units. The roots of an element of $G$ are the eigenvalues of its adjoint. An algebraic unit is a non-zero complex number that is an algebraic integer over $\mathbb{Z}$ as well as its inverse.
\end{definition}
If $G$ is solvable, G. Meigniez show the following theorem
\begin{theorem}{\cite{7}}
Let $G$ be a solvable Lie group and $\Gamma$ a subgroup of $G$ of finite type. If $\Gamma$ contains a subgroup that is polycyclic and uniform at the same time, then $\Gamma$ is a holonomy group of $G$-foliation.
\end{theorem}
In the completely solvable case we prove the following theorem
\begin{theorem}
Let $G$ be a completely solvable Lie group. $\Gamma$ a subgroup of $G$ of finite type and uniform. Then $\Gamma$ is polycyclic and uniform if and only if $\Gamma$ is virtually a holonomy group of a homogenius $G$-foliation.
\end{theorem}
Before giving the proof we have a property of rigidity of $G$
\begin{lemma}
Are $G_1$ and $G_2$ two simply connected Lie groups and completely solvable. Let $\Gamma$ be a discrete uniform subgroup in $G_1$. So everything
continuous homomorphism: $\rho:\Gamma\longrightarrow G_2$ extends in a single way into a homomorphism continue: $\widetilde{\rho}:G_1\longrightarrow G_2$.
\end{lemma}
Now we give the proof of the theorem
\begin{proof}
Let $G$ be a completely solvable Lie group and $\Gamma$ a subgroup of $G$.\\
And suppose that $\Gamma$ is of finite type and uniform. Every completely solvable Lie group is solvable Lie group, now we can use the Meigniez theorem \cite{7} and so $\Gamma$ is virtually a holonomy group of a homogenius $G$-foliation (because the Meigniez theorem hold on solvable case).\\
Conversely, assuming $\Gamma$ is polycyclic and uniform. Let $\varphi: \Gamma_H \longrightarrow\Gamma'\subset G$ where $\Gamma'$ is a subgroup of $\Gamma$ of finite index and $\Gamma_H$ a cocompact lattice of completely solvable Lie group $H$. Since $G$ is completely solvable and $\Gamma'$ is uniform (because lattice in completely solvable group is discret uniform subgroup), the continuous homomorphism: $\widetilde{\varphi}$  extends in a single way into a homomorphism continue: $\widetilde{\widetilde{\varphi}}:H\longrightarrow G$.\\
Then $ker\widetilde{\widetilde{\varphi}}$ is a homogenius Lie $G$-foliation on the manifold $H_{/{\Gamma_H}}$. And its holonomy group is $\widetilde{\widetilde{\varphi}}(\Gamma_H)=\Gamma'$
\end{proof}
\begin{theorem}
Let $G$ be simply connected and completely solvable Lie group, $\mathcal{F}$ a Lie $G$-foliation on a compact manifold $M$ and $\Gamma$ its holonomy group.\\
If $\Gamma$ is dicret then $\mathcal{F}$ is reverse image of homogenius Lie foliation.
\end{theorem}
we give this lemma before the proof of the theorem.
\begin{lemma}
Let $G$ be a simply connected and completely solvable Lie group and $\Gamma$ a closed subgroup of $G$. Then there exists a unique subgroup of Lie $S$ closed and connected of $G$ such that $\Gamma\in S$ and $ S/_{\Gamma}$ is compact.
\end{lemma}
Proof of the previous theorem.
\begin{proof}
Let $\mathcal{F}$ be Lie $G$-foliation on a compact manifold $M$ and $\Gamma$ its holonomy group.\\
Suppose that $\Gamma$ a discret subgroup of $G$ its holonomie group.\\
The holonomy group $\Gamma$ of such foliation $\mathcal{F}$ is a uniform subgroup of the finite type of the completely solvable Lie group, simply connected $G$ defining the transverse structure of $\mathcal{F}$.\\
Using the previous lemma there exists a unique subgroup of Lie $S$ closed and connected of $G$ such that $\Gamma\in S$ and $S/_{\Gamma}$ is compact. So $\Gamma$ is a lattice in $S$ (because lattice in completely solvable Lie group is discret and uniform subgroup).\\
The canonical injection $i:\Gamma\longrightarrow G$ extends in a single way into a surjective homomorphism such that the following diagram is commutative

\[
\xymatrix{
\Gamma \ar[r]^{j} \ar[dr]_{i} & S \ar[d]^{h} \\
 & G
}
\]
$$i=h\circ j$$
We have the shot exact sequence $1\longrightarrow F\longrightarrow S\longrightarrow G\longrightarrow 1$ which define a homogenius and classifiant Lie foliation $\mathcal{F}'$ on the manifold $S_{/\Gamma}$.\\
The holonomy group of $\mathcal{F}'$ is $h(\Gamma)=h\circ j(\Gamma)\approx \Gamma$. Then $\mathcal{F}$ is reverse image of $\mathcal{F}'$\cite{6}

\end{proof}

\end{document}